\documentclass[reqno]{amsart}

\usepackage{amssymb,latexsym}
\usepackage{amsfonts}
\usepackage{dsfont}
\usepackage[pdftex]{hyperref}
\usepackage{color}
\usepackage{graphicx}
\usepackage{parskip}
\usepackage{float}

\theoremstyle{plain}
\newtheorem{thm}{Theorem}

\numberwithin{equation}{section}

\newcommand{\p}{\partial}
\newcommand{\m}{\mathbb}

\newcommand{\g}{\gamma}
\newcommand{\ph}{\varphi}

\begin{document}

\title[k-abc]{Decay Properties of Solutions to a \\  4-parameter family of wave equations}
\author{Ryan C. Thompson}
\date{October 5, 2016}
\address{Department of Mathematics\\University of North Georgia\\ Dahlonega, GA 30597}
\email{ryan.thompson@ung.edu}



\begin{abstract}
In this paper, persistence properties of solutions are investigated for a 4-parameter family ($k-abc$ equation) of evolution equations having
$(k+1)$-degree nonlinearities  and containing
 as its integrable members the Camassa-Holm, the Degasperis-Procesi, Novikov and Fokas-Olver-Rosenau-Qiao equations.
 These properties will imply that strong solutions of the $k-abc$ equation will decay at infinity in the spatial variable provided that the initial data does.  Furthermore, 
it is shown that the equation exhibits unique continuation
for appropriate values of the parameters $k$, $a$, $b$, and $c$.
\end{abstract}

\maketitle

\section{Introduction}
For $k \in \m{Z}^+$, $k \geq 2$ and $a, b, c \in \m{R}$, we consider the Cauchy problem for the following $k-abc$ family of equations
\begin{align}
\label{kabc}
u_t &+u^ku_x-au^{k-2}u_x^3+D^{-2}\p_x\left[\frac{b}{k+1}u^{k+1}+cu^{k-1}u_x^2-a(k-2)u^{k-3}u_x^4\right] \nonumber \\ 
& \ +D^{-2}\left[[k(k+2)-8a-b-c(k+1)]u^{k-2}u_x^3-3a(k-2)u^{k-3}u_x^3u_{xx}\right]=0
\end{align}•
where $D^{-2} \doteq (1-\p_x^2)^{-1}$, and study its persistence properties.  We show that if the initial data is endowed with exponential decay at infinity, then the corresponding solution will carry this property.
When $a=0$ and $k$ a positive odd integer with $k\geq 1$ and $b \in [0,k(k+2)]$ then,
 utilizing the aforementioned behavior of strong solutions, we show that this equation exhibits unique continuation in the sense that if the initial value $u(x,0)$ is given the property of decaying exponentially, then the solution $u(x,t)$ must be identically zero if assumed to be decaying exponentially at some later time $t>0$.   These are natural extensions of the results 
 proved in Himonas, Misio\l ek, Ponce and Zhou \cite{hmpz} for the Camassa-Holm equation.

The $k-abc$ equation was first studied by Himonas and Mantzavinos \cite{hm1} where they showed well-posedness in Sobolev spaces $H^s$ for $s>5/2$.  They also provided a sharpness result on the data-to-solution map and proved that it is not uniformly continuous from any bounded subset of $H^s$ into $C([0,T];H^s)$.  It was shown, however, that the solution map is H\"{o}lder continuous if $H^s$ is equipped with a weaker $H^r$ norm where $r \in [0,s)$.  The equation was also studied in Barostichi, Himonas and Petronilho in \cite{bhp} where they exhibited a power series method in abstract Banach spaces provided analytic initial data, thereby establishing a Cauchy-Kovalevsky type theorem for the $k-abc$ equation \eqref{kabc}.

It is important to note that when Himonas and Mantzavinos ushered in the $k-abc$ equation \eqref{kabc}, they unified many Camass-Holm type equations.  For instance, when $k=2$ and $c = (6-6a-b)/2$ you obtain the $ab$-equation
\begin{equation}
\label{ab}
u_t+u^2u_x-au_x^3+D^{-2}\p_x\left[\frac b3u^3+\frac{6-6a-b}{2}uu_x^2\right]+D^{-2}\left[\frac{2a+b-2}{2}u_x^3\right]=0,
\end{equation}
which was also introduced by Himonas and Mantzavinos in \cite{hm2} and was shown to possess both periodic and non-periodic traveling wave solutions.  In fact, the non-periodic multipeakon traveling wave solutions to the $ab$-equation \eqref{ab} take the form
\[
u(x,t) = \sum_{j=1}^np_j(t)e^{-|x-q_j(t)|},
\]
provided that the positions $q_j$ and momenta $p_j$ satisfy an appropriate system of nonlinear differential equations.

When $a=0$ and $c=(3k-b)/2$, the $k-abc$ equation yields the generalized Camassa-Holm equation (g-$kb$CH)
\begin{equation}
\label{gkbch}
u_t+u^ku_x+D^{-2}\p_x\left[\frac {b}{k+1}u^{k+1}+\frac{3k-b}{2}u^{k-1}u_x^2\right]+D^{-2}\left[\frac{(k-1)(b-k)}{2}u^{k-2}u_x^3\right]=0.
\end{equation}
Well-posedness in $H^s$ with $s>3/2$ was shown in \cite{hh1}, while its multipeakon traveling waves were derived in \cite{gh}.  Furthermore, Himonas and Thompson \cite{ht} explored asymptotic behavior of strong solutions, unique continuation and conditions needed to admit global solutions.

If we let $k=1$ in the aforementioned then we also obtain the $b$-family of equations with quadratic nonlinearities which was introduced by Holm and Staley in \cite{hs1} and \cite{hs2}.  We may also see its local form as
\begin{equation}
\label{bfam}
\underbrace{m_t}_\text{evolution}+\underbrace{um_x}_\text{convection}+\underbrace{bu_xm}_\text{stretching}=0, \ \ \ m = u-u_{xx},
\end{equation}
which expresses a fine balance between, evolution, convection and stretching.  It is within this balance that one may find a starting point for not only deriving peakon solutions, but also establishing infinite speed of propagation and other asymptotic behavior of solutions.

When we let $a=1/3$ and $b=2$ in the $ab$-equation \eqref{ab} we obtain the Fokas-Olver-Rosenau-Qiao (FORQ) equation
\begin{equation}
\label{forq}
u_t+u^2u_x-\frac13u_x^3+D^{-2}\p_x\left[\frac23u^3+uu_x^2\right]+D^{-2}\left[\frac13u_x^3\right]=0,
\end{equation}
which was derived in a variety of ways by Fokas \cite{fokas}, Olver and Rosenau \cite{or} and Qiao \cite{qiao} and also appeared in a work by Fuchsteiner \cite{fuch}.

Observe that for  $k=1$ and $b=2$ the g-$kb$CH equation \eqref{gkbch}
gives the celebrated Camassa-Holm (CH) equation
\begin{equation}
\label{ch}
(1-\p_x^2)u_t=uu_{xxx}+2u_xu_{xx}-3uu_x
\end{equation}
which appears in the context of hereditary symmetries studied by Fokas and Fuchssteiner \cite{ff}.  This equation was first written explicitly and derived from the Euler equations  by Camassa and Holm in \cite{cam}, where they also found its 
``peakon" traveling wave solutions.
These are solutions with a discontinuity in the first spacial derivative at its crest. The simplest one in the non-periodic case
is of the form $u_c(x, t)=ce^{-|x-t|}$, where $c$ is a positive constant.
The CH equation possesses many  other 
remarkable properties such as infinitely many conserved quantities, a bi-Hamiltonian structure and a Lax pair.
For more information about how CH arises in the context of hereditary symmetries we refer to \cite{ff}.  Concerning it's physical relevance, we refer the reader to the works by Johnson \cite{j1}, \cite{j2} and Constantine and Lannes \cite{const1}.

For $k=1$ and $b=3$ in \eqref{gkbch}, we obtain the Degasperi-Procesi (DP) equation
\begin{equation}
\label{dp}
(1-\p_x^2)u_t=uu_{xxx}+3u_xu_{xx}-4uu_x.
\end{equation}
Degasperi and Procesi \cite{dp} discovered this equation in 1998 as one of the three equations to satisfy asymptotic integrability conditions in the following family of equations
\begin{equation}
\label{fam}
u_t+c_0u_x+\g u_{xxx}-\alpha^2u_{txx}=(c_{1,n}u^2+c_{2,n}u_x^2+c_3uu_{xx})_x,
\end{equation}
where $\alpha, c_0, c_{1,n}, c_{2,n}, c_3 \in \m{R}$ are constants.  Other integrable members of the family \eqref{fam}
are the CH and the KdV equations.  Furthermore, the DP and CH are the only integrable members of the $b$-family equation which is obtained from the g-$kb$CH equation by letting $k=1$ and $b \in \m{R}$ (see \cite{mn}).  Also, like the CH equation \eqref{ch}, the DP equation \eqref{dp} possesses peakon solutions of the form $u_c(x,t)=ce^{-|x-ct|}$.

Recently, Vladimir Novikov \cite{nov}, in addition to 
several integrable equations with quadratic nonlinearities
like the CH and DP, generated 
about 10 integrable equations with cubic nonlinearities while investigating the integrability of Camassa-Holm type equations of the form 
\begin{equation}
\label{gen}
(1-\p_x^2)u_t=P(u,u_x,u_{xx},...),
\end{equation}
where $P$ is a polynomial of $u$ and its derivatives.  
One of these CH-type equations  with cubic nonlinearities, 
which happened to be a new integrable equation, is  the following
\begin{equation}
\label{ne}
(1-\p_x^2)u_t=u^2u_{xxx}+3uu_xu_{xx}-4u^2u_x,
\end{equation}
which can be found from the g-$kb$CH equation \eqref{gkbch} by letting $k=2$ and $b=3$.  Equation \eqref{ne} is now called the Novikov equation (NE) and also possesses many properties exhibited by Camassa-Holm type equations.  One of which is the existence of peakon solutions of the form $u_c(x,t)=\sqrt{c}e^{-|x-ct|}$.  Unlike CH and DP, however, unique continuation has not been shown for NE due to difficulties in handling its cubic nonlinearities.

In the direction of peakon traveling wave solutions, 
it has been shown in  \cite{hm2}  that the $k-abc$ equation \eqref{kabc} possesses these solitary wave solutions on the circle and the real line.  On the real line, these solutions take the form
\[
u(x,t)=\g e^{-|x-(1-a)\g^kt|}, \ \ \g \in \m{R},
\]
where $(1-a)\g^k$ is the wave speed.  On the circle, these solutions take the form
\[
u(x,t)=\g\cosh\left(\left[x-\left[1+(1-a)\sinh^2\pi\right]\cosh^{k-2}(\pi)\g^kt\right]_p-\pi\right), \ \ \g \in \m{R},
\]
where
\[
[x]_p \doteq x-2\pi\left\lfloor \frac{x}{2\pi} \right\rfloor
\]
provided that $6a+b+2c=3k$.  Furthermore, under appropriate conditions on the four parameters the $k-abc$ equation also admits multi-peakon solutions.

In regards to conservation laws, we take note that the $H^1$ norm is conserved under particular parameter conditons.  In fact, we see have
\[
\frac12\frac{d}{dt}\|u(t)\|_{H^1}=\int\left(\left[2k-2c-1-\frac2k(9a+b+2c-3k)\right]u^ku_xu_{xx}-\frac{3a}{4}(k-1)(k-2)u^{k-3}u_x^5\right)dx
\]
and so we find the $H^1$ norm is conserved if $9a+b+4c=9$ when $k=2$ and if $a=0$ and $2c+\frac2k(b+2c-3k)+1=2k$ when $k \geq 3$.

The well-posedness of the $k-abc$ equation in Sobolev spaces and the continuity properties of its solution map have been studied in \cite{hm1}.  There, the following result was proved.
\begin{thm}
\label{t1}
Let $a,b,c \in \m{R}$ with $a \neq 0$ and $k \in \m{N}$ with $k \geq 2$.  Then the Cauchy problem for the $k-abc$ equation \eqref{kabc} with initial condition $u(x,0)=u_0(x) \in H^s$. $x \in \m{R}$ or $\m{T}$, $s > 5/2$, has a unique solution $u \in C([0,T];H^s)$ which admits the estimate
\begin{equation}
\|u(t\|_{H^s} \leq 2^{1+\frac 1k}\|u_0\|_{H^s}, \ \ 0 \leq t \leq T \leq (2^{k+1}kc_s\|u_0\|_{H^s}^k)^{-1},
\end{equation}•
where the constant $c_s>0$ depends only on $s$.  Furthermore, the data-to-solution map $u_0 \mapsto u(t)$ is continuous.
\end{thm}

It is also important to note that when $a=0$, the $k-abc$ equation closely resembles the g-$kb$CH equation \eqref{gkbch} and a similar well-posedness result was proven, albeit for $s>3/2$, in \cite{hh1}.

From the above well-posedness result, we may now extract our persistence properties and unique continuation.  Our
basic assumption is that the initial data as well as its first spacial derivative  decay exponentially.  We also endow the solution at time $t>0$ with an analogous condition which guarantees that our solution will be identically zero.  The following results are based on the work
in  \cite{hmpz} for the CH equation.
\begin{thm}
\label{t2}
Assume that $T>0$, $s>5/2$ and $u \in C([0,T];H^s)$ is a strong solution of \eqref{kabc}.  If the initial data satisfies certain decay conditions at infinity, more precisely, if there is some $\theta \in (0,1)$ such that
\[
|u_0|, \ |\p_xu_0| \sim O(e^{-\theta |x|}) \ \ as \ \ |x| \to \infty,
\]
then 
\[
|u(x,t)|, \ |\p_xu(x,t)| \sim O(e^{-\theta |x|}) \ \ as \ \ |x| \to \infty
\]
uniformly with respect to $t \in [0,T]$.
\end{thm}

\emph{Notation.}  We shall say that
\[
|f(x)| \sim O(e^{\alpha x}) \ \ \text{as} \ \ x \to \infty \ \ \text{if} \ \ \lim_{x \to \infty}\frac{|f(x)|}{e^{\alpha x}}=L,
\]
for some $L>0$, and
\[
|f(x)| \sim o(e^{\alpha x}) \ \ \text{as} \ \ x \to \infty \ \ \text{if} \ \ \lim_{x \to \infty}\frac{|f(x)|}{e^{\alpha x}}=0.
\]

\begin{thm}
\label{kabc-uniq-cont}
Assume that for some $T>0$ and $s>5/2$,
$
u \in C([0,T];H^s(\m{R}))
$
is a solution to the $k-abc$ initial value problem, $k$ is an odd integer with $k \geq 1$ and that 
\[
a=0, \ \  c=\frac{k(k+2)-b}{k+1}, \ \ b \in [0, k(k+2)].
\]
If $u_0(x)=u(x,0)$ satisfies that for some $\alpha \in \left(\frac{1}{k+1},1\right)$,
\begin{equation}
\label{decay1}
|u_0(x)| \sim o(e^{-x}), \ \ \ \ and \ \ \ \ |\p_xu_0(x)| \sim O(e^{-\alpha x}) \ \ \text{as} \ x \rightarrow \infty,
\end{equation}
and there exists $t_1 \in (0,T]$ such that
\begin{equation}
\label{decay2}
|u(x,t_1)| \sim o(e^{-x}) \ \ as \ \ x \rightarrow \infty,
\end{equation}
then $u \equiv 0$.
\end{thm}

The paper is organized as follows.  In Section 2, we prove the persistence properties of the $k-abc$ equation as listed in Theorem \ref{t2}.  Then in Section 3 we use Theorem \ref{t2} to prove Theorem \ref{kabc-uniq-cont}.  Finally, we make a remark about the aymptotic profile of solutions when provided compactly supported initial data.

%
%
%
%
\section{Proof of Theorem \ref{t2}}

\begin{proof}
Assume that $u \in C([0,T];H^s)$ is a strong solution to \eqref{kabc} with $s > 5/2$.  Let
\begin{align}
\label{fs}
&F_1(u)=\frac{b}{k+1}u^{k+1}+cu^{k-1}u_x^2-a(k-2)u^{k-3}u_x^4, \\ &F_2(u)=[k(k+2)-8a-b-c(k+1)]u^{k-2}u_x^3-3a(k-2)u^{k-3}u_x^3u_{xx}.
\end{align}•
Then we have that the $k-abc$ equation \eqref{kabc} takes the form 
\begin{equation}
\label{nonloc2}
u_t+u^ku_x-au^{k-2}u_x^3+\p_xG*F_1(u)+G*F_2(u)=0.
\end{equation}•
Also, let
\begin{equation}
\label{m}
M=\sup_{t \in [0,T]}\|u(t)\|_{H^s},
\end{equation}•
\begin{equation}
\label{phi}
\ph_N(x)=
\begin{cases}
e^{\theta|x|}, \ \ |x|<N \\ 
e^{\theta N}, \ \ \  |x| \geq N,
\end{cases}•
\end{equation}•
where $N \in \m{N}$ and $\theta \in (0,1)$.  Observe that for all $N$ we have
\[
0 \leq |\ph_N'(x)| \leq \ph_N, \ \ \ a.e. \ \ x \in \m{R}.
\]
Multiplying \eqref{nonloc2} by $(u\ph_N)^{2p-1}\ph_n$ for $p \in \m{Z}^+$ and integrating over the real line we obtain
\begin{align}
\label{24}
\frac{1}{2p}\frac{d}{dt}&\int_\m{R}(u\ph_N)^{2p}dx+\int_\m{R}u^ku_x(u\ph_N)^{2p-1}\ph_Ndx \nonumber \\
&-\int_\m{R} au^{k-2}u_x^3(u\ph_N)^{2p-1}\ph_N dx \nonumber \\ 
&+\int_\m{R}(\p_xG*F_1(u))(u\ph_N)^{2p-1}\ph_Ndx \nonumber \\ 
&+\int_\m{R}(G*F_2(u))(u\ph_N)^{2p-1}\ph_Ndx=0.
\end{align}•
By the Sobolev embedding theorem and \eqref{m} we have
\[
\|u\|_{L^\infty}+\|u_x\|_{L^\infty}+\|u_{xx}\|_{L^\infty} \leq CM.
\]
This leads us to achieve the following estimates
\begin{align}
\label{25}
&\left|\int_\m{R}u^ku_x(u\ph_N)^{2p-1}\ph_Ndx\right| \leq \|u\|_{\infty}^{k-1}\|u_x\|_{\infty}\|u\ph_N\|_{2p}^{2p} \leq CM^k\|u\ph_N\|_{2p}^{2p} \\ 
&\left|a\int_\m{R} u^{k-2}u_x^3(u\ph_N)^{2p-1}\ph_N dx\right| \leq |a|\|u\|_{\infty}^{k-3}\|u_x\|_{\infty}^3\|u\ph_N\|_{2p}^{2p} \leq CM^k\|u\ph_N\|_{2p}^{2p}.
\end{align}•
By H\"{o}lder's inequality, we have the following estimates
\begin{align}
\label{26}
&\left|\int_\m{R}(\p_xG*F_1(u))(u\ph_N)^{2p-1}\ph_Ndx\right| \leq \|u\ph_N\|_{2p}^{2p-1}\|(\p_xG*F_1(u))\ph_N\|_{2p}, \nonumber \\ 
&\left|\int_\m{R}(G*F_2(u))(u\ph_N)^{2p-1}\ph_Ndx\right| \leq \|u\ph_N\|_{2p}^{2p-1}\|(G*F_2(u))\ph_N\|_{2p}.
\end{align}•
From \eqref{24} and the above estimates, this implies
\begin{equation}
\label{ode1}
\frac{d}{dt}\|u\ph_N\|_{2p} \leq CM^k\|u\ph_N\|_{2p}+\|(\p_xG*F_1(u))\ph_N\|_{2p}+\|(G*F_2(u))\ph_N\|_{2p}.
\end{equation}•
By Gronwall's inequality, \eqref{ode1} implies the following estimate
\begin{equation}
\label{est1}
\|u\ph_N\|_{2p} \leq \left(\|u_0\ph_N\|_{2p}+\int_0^t\left[\|(\p_xG*F_1(u))\ph_N\|_{2p}+\|(G*F_2(u))\ph_N\|_{2p}\right]d\tau\right)e^{CM^kt}.
\end{equation}•
Now differentiating \eqref{nonloc2} with respect to the spacial variable $x$, multiplying by $(\ph_Nu_x)^{2p-1}\ph_N$ and integrating over the real line yields
\begin{align}
\label{29}
\frac{1}{2p}\frac{d}{dt}&\int_\m{R}(u_x\ph_N)^{2p}dx+\int_\m{R}ku^{k-1}(u_x)^2(\ph_Nu_x)^{2p-1}\ph_Ndx \nonumber \\ 
&+\int_\m{R}u^ku_{xx}(\ph_Nu_x)^{2p-1}\ph_Ndx-\int_\m{R}a(k-2)u^{k-3}u_x^4(\ph_Nu_x)^{2p-1}\ph_N dx \nonumber \\ 
&-\int_\m{R}3au^{k-2}u_x^2u_{xx}(\ph_Nu_x)^{2p-1}\ph_N dx+\int_\m{R}(\p_x^2G*F_1(u))(\ph_Nu_x)^{2p-1}\ph_Ndx \nonumber \\ 
&+\int_\m{R}(\p_xG*F_2(u))(\ph_Nu_x)^{2p-1}\ph_Ndx=0.
\end{align}•
This leads us to obtain the following estimates
\begin{align}
\label{210}
&\left|\int_\m{R}ku^{k-1}(u_x)^2(\ph_Nu_x)^{2p-1}\ph_Ndx\right| \leq CM^k\|u_x\ph_N\|_{2p}^{2p}, \nonumber \\
&\left|\int_\m{R}a(k-2)u^{k-3}u_x^4(\ph_Nu_x)^{2p-1}\ph_N dx\right| \leq  CM^k\|u_x\ph_N\|_{2p}^{2p} \nonumber \\
&\left|\int_\m{R}3au^{k-2}u_x^2u_{xx}(\ph_Nu_x)^{2p-1}\ph_N dx\right| \leq CM^k\|u_x\ph_N\|_{2p}^{2p} \nonumber \\
&\left|\int_\m{R}(\p_x^2G*F_1(u))(\ph_Nu_x)^{2p-1}\ph_Ndx\right| \leq \|u_x\ph_N\|_{2p}^{2p-1}\|(\p_x^2G*F_1(u))\ph_N\|_{2p}, \nonumber \\ 
&\left|\int_\m{R}(\p_xG*F_2(u))(\ph_Nu_x)^{2p-1}\ph_Ndx\right| \leq \|u_x\ph_N\|_{2p}^{2p-1}\|(\p_xG*F_2(u))\ph_N\|_{2p}.
\end{align}•
For the third integral in \eqref{29}, we estimate as follows
\begin{align}
\label{211}
\left|\int_\m{R}u^ku_{xx}(\ph_Nu_x)^{2p-1}\ph_Ndx\right|&=\left|\int_\m{R}u^k(\ph_N u_x)^{2p-1}(\p_x(\ph_Nu_x)-\ph_N' u_x)dx\right| \nonumber \\ 
&=\left|\int_\m{R}u^k(\ph_Nu_x)^{2p-1}\p_x(\ph_Nu_x)dx-\int_\m{R}u^k(\ph_Nu_x)^{2p-1}\ph_N'u_xdx\right| \nonumber \\ 
&=\left|\frac{1}{2p}\int_\m{R}u^k\p_x[(\ph_Nu_x)^{2p}]dx-\int_\m{R}u^k(\ph_Nu_x)^{2p-1}\ph_N'u_xdx\right| \nonumber \\ 
&=\left|\frac{-k}{2p}\int_\m{R}u^{k-1}u_x(\ph_Nu_x)^{2p}dx-\int_\m{R}u^k(\ph_Nu_x)^{2p-1}\ph_N'u_xdx\right| \nonumber \\ 
& \leq CM^k\|u_x\ph_N\|_{2p}^{2p}.
\end{align}•
From \eqref{29} - \eqref{211}, we achieve the following differential inequality
\begin{equation}
\label{ode2}
\frac{d}{dt}\|u_x\ph_N\|_{2p} \leq CM^k\|u_x\ph_N\|_{2p}+\|(\p_x^2G*F_1(u))\ph_N\|_{2p}+\|(\p_xG*F_2(u))\ph_N\|_{2p}.
\end{equation}•
From \eqref{ode2} and Gronwall's inequality, we obtain
\begin{equation}
\label{est2}
\|u_x\ph_N\|_{2p} \leq \left(\|\p_xu_0\ph_N\|_{2p}+\int_0^t[\|(\p_x^2G*F_1(u))\ph_N\|_{2p}+\|(\p_xG*F_2(u))\ph_N\|_{2p}]d\tau\right)e^{CM^kt}
\end{equation}•
By adding \eqref{est1} and \eqref{est2}, we have the following
\begin{align}
\label{est3}
\|u\ph_N\|_{2p}+\|u_x\ph_N\|_{2p} \ \leq & \ (\|u_0\ph_N\|_{2p}+\|\p_xu_0\ph_N\|_{2p})e^{CM^kt} \nonumber \\ 
&+\left(\int_0^t[\|(\p_xG*F_1(u))\ph_N\|_{2p}+\|(G*F_2(u))\ph_N\|_{2p}]d\tau\right)e^{CM^kt} \nonumber \\ 
&+\left(\int_0^t[\|(\p_x^2G*F_1(u))\ph_N\|_{2p}+\|(\p_xG*F_2(u))\ph_N\|_{2p}]d\tau\right)e^{CM^kt}
\end{align}•

Now, for any function $f \in L^1 \cap L^\infty$, $\lim_{q \to \infty}\|f\|_q=\|f\|_\infty$.  Since we have that $F_i(u) \in L^1 \cap L^\infty$ (for $i=1,2$) and $G \in W^{1,1}$, we know that $\p_x^jG*F_i(u) \in L^1 \cap L^\infty$ (for $i=1,2$ and $j=0,1,2$).  Thus, by taking the limit of \eqref{est3} as $p \to \infty$, we obtain
\begin{align}
\label{est4}
\|u\ph_N\|_{\infty}+\|u_x\ph_N\|_{\infty} \ \leq & \ (\|u_0\ph_N\|_{\infty}+\|\p_xu_0\ph_N\|_{\infty})e^{CM^kt} \nonumber \\ 
&+\left(\int_0^t[\|(\p_xG*F_1(u))\ph_N\|_{\infty}+\|(G*F_2(u))\ph_N\|_{\infty}]d\tau\right)e^{CM^kt} \nonumber \\ 
&+\left(\int_0^t[\|(\p_x^2G*F_1(u))\ph_N\|_{\infty}+\|(\p_xG*F_2(u))\ph_N\|_{\infty}]d\tau\right)e^{CM^kt}
\end{align}•

A simple calculation, which we give in the appendix, shows that for $\theta \in (0,1)$
\begin{equation}
\label{calc1}
\ph_N(x)\int_\m{R}e^{-|x-y|}\frac{1}{\ph_N(y)}dy \leq \frac{4}{1-\theta}=C_0.
\end{equation}•
Thus, for any functions $f,g \in L^\infty$ we have
\begin{align}
\label{calc2}
|\ph_NG*(f^2g)(x)|&=\left|\frac12\ph_N\int_\m{R}e^{-|x-y|}f^2(y)g(y)dy\right| \nonumber \\ 
&\leq \frac12\left(\ph_N\int_\m{R}e^{-|x-y|}\frac{1}{\ph_N(y)}dy\right)\|f\|_\infty\|g\|_\infty\|f\ph_N\|_\infty \nonumber \\ 
&\leq C_0\|f\|_\infty\|g\|_\infty\|f\ph_N\|_\infty.
\end{align}•
Similarly, we have
\begin{equation}
\label{calc3}
|\ph_N\p_xG*(f^2g)(x)| \leq C_0\|f\|_\infty\|g\|_\infty\|f\ph_N\|_\infty
\end{equation}•
and
\begin{equation}
\label{calc4}
|\ph_N\p_x^2G*(f^2g)(x)|=|\ph_N[G*(f^2g)-f^2g]| \leq C_0\|f\|_\infty\|g\|_\infty\|f\ph_N\|_\infty.
\end{equation}•
Therefore, with $f=u_x$ we obtain the following estimates
\begin{align}
\label{est5}
&\|(G*F_2(u))\ph_N\|_{\infty} \leq CM^k\|u_x\ph_N\|_\infty, \nonumber \\ 
&\|(\p_xG*F_2(u))\ph_N\|_{\infty} \leq CM^k\|u_x\ph_N\|_\infty.
\end{align}•
Then, by following the same procedure as in \eqref{calc2} we obtain
\begin{align}
\label{est6}
&\|(\p_xG*F_1(u))\ph_N\|_{\infty} \leq CM^k\|u\ph_N\|_\infty, \nonumber \\ 
&\|(\p_x^2G*F_1(u))\ph_N\|_{\infty} \leq CM^k\|u\ph_N\|_\infty.
\end{align}•
So, by estimates \eqref{est4}, \eqref{est5} and \eqref{est6} we achieve the following
\begin{align}
\|u\ph_N\|_{\infty}+\|u_x\ph_N\|_{\infty} &\leq C(\|u_0\ph_N\|_{\infty}+\|\p_xu_0\ph_N\|_{\infty}) \nonumber \\ 
&+C\int_0^t(\|u\ph_N\|_\infty+\|u_x\ph_N\|_\infty)d\tau.
\end{align}•
Hence, for any $N \in \m{Z}^+$ and any $t \in [0,T]$, we have by Gronwall's inequality that
\begin{align}
\|u\ph_N\|_{\infty}+\|u_x\ph_N\|_{\infty} &\leq  C(\|u_0\ph_N\|_{\infty}+\|\p_xu_0\ph_N\|_{\infty}) \nonumber \\ 
&\leq  C(\|u_0e^{\theta|x|}\|_{\infty}+\|\p_xu_0e^{\theta|x|}\|_{\infty}).
\end{align}
By taking the limit as $N \to \infty$ we have that for any $t \in [0,T]$
\begin{equation}
|u(x,t)|e^{\theta|x|}+|u_x(x,t)|e^{\theta|x|} \leq  C(\|u_0e^{\theta|x|}\|_{\infty}+\|\p_xu_0e^{\theta|x|}\|_{\infty}).
\end{equation}
This concludes our proof of Theorem \ref{t2}.
\end{proof}

\textbf{Remark 2.1:}  In fact, let $\theta \in (0,1)$ and $j=0,1,2$, if the initial data $u(x,0) = u_0$ satisfy
\[
\p_x^ju_0 \sim O(e^{-\theta|x|}), \ \ as \ \ x \to \infty,
\]
then the solution $u(x,t)$ carries the same exponential decay properties, i.e.
\[
\p_x^ju \sim O(e^{-\theta|x|}), \ \ as \ \ x \to \infty.
\]
Theorem \ref{t2} tells us that the solution $u(x,t)$ can decay as $e^{-\theta|x|}$ as $x \to \infty$ for $\theta \in (0,1)$.  The next theorem provides us with further decay properties on the solution if the only assumptions we make on $k$, $a$, $b$, and $c$ are those stated in the introduction.

%
%
%
%
\section{Proof of Theorem \ref{kabc-uniq-cont}}

First, we will assume the result from Theorem \ref{t2}.  Integrating \eqref{nonloc2} over the time interval $[0,t_1]$, noting that $(1-\p_x^2)^{-1}f=G*f$ for all $f \in L^2$ where $G=\frac12e^{-|x|}$ and assuming the stated restrictions on parameters $k$, $a$, $b$ and $c$ we get
\begin{equation}
\label{9911}
u(x,t_1)-u(x,0)+\int_0^{t_1}u^k\p_xu(x,\tau)d\tau+\int_0^{t_1}\p_xG*F_1(u)(x,\tau)d\tau=0,
\end{equation}
where
\[
F_1(u)=\frac{b}{k+1}u^{k+1}+cu^{k-1}u_x^2.
\]
By hypothesis \eqref{decay1} and \eqref{decay2} we have
\begin{equation}
\label{9912}
u(x,t_1)-u(x,0) \sim o(e^{-x}) \ \ \text{as} \ \ x \rightarrow \infty.
\end{equation}
From \eqref{decay1} and Theorem \ref{t2} it follows that
\begin{equation}
\label{9913}
\int_0^{t_1}u^k\p_xu(x,\tau)d\tau \sim O(e^{-\alpha(k+1) x}) \ \ \text{as} \ \ x \rightarrow \infty,
\end{equation}
and so
\begin{equation}
\label{9914}
\int_0^{t_1}u^k\p_xu(x,\tau)d\tau \sim o(e^{-x}) \ \ \text{as} \ \ x \rightarrow \infty.
\end{equation}
We shall show that if $u \neq 0$, then the last two terms in \eqref{9911} is $O(e^{-x})$ but not $o(e^{-x})$.  Thus, we have
\begin{align}
\label{9915}
\int_0^{t_1}\p_xG*F_1(u)(x,\tau)d\tau=\p_xG*\int_0^{t_1}F_1(u)(x,\tau)d\tau=\p_xG*f_1(x),
\end{align}
where by \eqref{decay1} and Theorem \ref{t2},
\begin{equation}
\label{9916}
0 \leq f_1(x) \sim O(e^{-\alpha(k+1) x}), \ \ \text{so that} \ \ f_1(x) \sim o(e^{-x}) \ \ \text{as} \ \ x \rightarrow \infty.
\end{equation}
Therefore,
\begin{equation}
\label{9917}
\p_xG*f_1(x)=-\frac12e^{-x}\int_{-\infty}^xe^yf_1(y)dy+\frac12e^x\int_x^\infty e^{-y}f_1(y)dy.
\end{equation}
From \eqref{9916} it follows that
\begin{equation}
\label{9918}
e^x\int_x^\infty e^{-y}f_1(y)dy=o(1)e^x\int_x^\infty e^{-2y}dy \sim o(1)e^{-x} \sim o(e^{-x}),
\end{equation}
and if $f_1 \neq 0$ one has that 
\begin{equation}
\label{9919}
\int_{-\infty}^xe^yf_1(y)dy \geq c_0, \ \ \ \text{for} \ \ x \gg 1, \ c_0>0.
\end{equation}
Hence, the first integral in \eqref{9917} satisfies
\begin{equation}
\label{9920}
\frac12e^{-x}\int_{-\infty}^xe^yf_1(y)dy \geq \frac{c_0}{2}e^{-x}, \ \ \text{for} \ \ x \gg 1,
\end{equation}
which combined with \eqref{9911}-\eqref{9913} yields a contradiction.  Thus, $f_1(x) \equiv 0$ and consequently $u \equiv 0$. 
$\Box$

\textbf{Remark on compactly supported initial data:}  In this section, we reflect on the property of unique continuation which we have just shown the Cauchy problem for the $k-abc$ equation \eqref{kabc} to exhibit.  In the case of compactly supported initial data
unique continuation is essentially infinite speed of propagation of its support.  Therefore, it is natural to ask the question:  How will strong solutions behave at infinity when given compactly supported initial data?  Our intentions are to follow the methods established in \cite{hmpz} and Himonas and Thompson \cite{ht}.   We typically need two ingredients in order to provide an asymptotic profile of our solution.

First, we observe that the $k-abc$ equation \eqref{kabc} may be written in the following local form
\begin{align}
\label{kabc_loc}
u_t &-u_{xxt}+(b+1)u^ku_x+(2c-3k)u^{k-1}u_xu_{xx}-u^ku_{xxx} \nonumber \\ 
&+(3k-9a-b-2c)u^{k-2}u_x^3+6au^{k-2}u_xu_{xx}^2+3au^{k-2}u_x^2u_{xxx}=0.
\end{align}

Our goal would now be to find the appropriate $m$-equation of the form
\begin{equation}
\label{kabc_m}
\underbrace{m_t}_\text{evolution}+\gamma_1\underbrace{u^km_x}_\text{convection}+\gamma_2\underbrace{(u^k)_xm}_\text{stretching}=0, \ \ \ m = u-u_{xx}, \ \ \gamma_1, \g_2 \in \m{R},
\end{equation}
which expresses a balance between evolution, convection and stretching.  Furthermore, we consider the following particle trajectories on our solution $u(x,t)$ given by
\begin{equation}
\label{parttraj}
\begin{cases}
\eta_t(x,t)=u^k(\eta(x,t),t) \\ 
\eta(x,0)=x,
\end{cases}
\end{equation}
and then seek to establish a following conservation law of the form
\begin{equation}
\label{cons_1}
m(\eta(x,t),t)\eta_x^{\g_3} = m_0(x), \ \ \g_3 \in \m{R}.
\end{equation}
With \eqref{kabc_m} to \eqref{cons_1}, we can establish an asymptotic profile of our solution provided compactly supported initial data.  To see this construction, we refer the reader to \cite{hmpz} and \cite{ht}.   In our search for the aforementioned with appropriate $\g_j$'s, however, we found that due to the coefficients of terms $u_t-u_{xxt}$, $(b+1)u^ku_x$ and $u^ku_{xxx}$ in our local equation \eqref{kabc_loc} forces our $\g_j$'s to be
\[
\g_1 = 1, \ \ \ \g_2 = \frac bk, \ \ \ \g_3 = \frac bk.
\]
Unfortunately, this forces $a=0$ and $3k-9a-b-2c=0$ which reduces us to the g-$kb$CH equation as studied in \cite{ht}.  So while it is possible to obtain the result on unique continuation for appropriate values of $a$, $b$ and $c$, the procurement of asymptotic profiles under compactly supported initial data is more delicate in regards to the previous literature.

\end{document}